\documentclass[12pt]{article}

\usepackage{amsmath,amssymb,amsfonts}

\newcommand{\Spec}{\mbox{Spec}\,}

\newcommand{\depth}{\mbox{depth}\,}
\renewcommand{\dim}{\mbox{dim}\,}

\newcommand{\Min}{\mbox{Min}}

\newcommand{\T}{\mathrm}

\newcommand{\fm}{\frak{m}}
\newcommand{\fp}{\frak{p}}
\newcommand{\fq}{\frak{q}}
\newcommand{\fn}{\frak{n}}

\begin{document}

\begin{center}

\LARGE{\bf Tensor Products of Some Special Rings\footnote{This
research was supported in part by a grant from IPM.}}

\end{center}

\begin{center}

{\bf Masoud Tousi$^{\it ab}$ and Siamak Yassemi$^{\it ac}$}
\end{center}

\begin{verse}
$(a)$ {\it Institute for Studies in Theoretical Physics and
Mathematics}

$(b)$ {\it Department of Mathematics, Shahid Beheshti University}

$(c)$ {\it Department of Mathematics, University of Tehran.}

\end{verse}

\vspace{.3in}

\begin{abstract}

In this paper we solve a problem, originally raised by
Grothendieck, on the properties, i.e. Complete intersection,
Gorenstein, Cohen--Macaulay, that are conserved under tensor
product of algebras over a field $k$.

\end{abstract}

\vspace{.2in}



\baselineskip=18pt

\vspace{.3in}

\section*{0. Introduction} Throughout this note all rings
and algebras considered in this paper are commutative with
identity elements, and all ring homomorphisms are unital.
Throughout, $k$ stands for a field.

Among local rings there is a well--known chain:

Regular $\Rightarrow$ Complete intersection $\Rightarrow$
Gorenstein $\Rightarrow$ Cohen--Macaulay.

These concepts are extended to non--local rings: for example a
ring is regular if for all prime ideal $\fp$ of $R$, $R_{\fp}$ is
a regular local ring.

In this paper, we shall investigate if these properties are
conserved under tensor product operations. It is well--known that
the tensor product $R\otimes_AS$ of regular rings is not regular
in general, even if we assume $R$ and $S$ are $A$--algebra and $A$
is a field, see Remark 1.7. In [{\bf 5}], Watanabe, Ishikawa,
Tachibana, and Otsuka, showed that under a suitable condition
tensor products of regular rings are complete intersections. It is
proved in [{\bf 3}], that the tensor product $R\otimes_AS$ of
Cohen--Macaulay rings are again Cohen--Macaulay if we assume $R$
is flat $A$--module and $S$ is a finitely generated $A$--module,
and in [{\bf 5}], it is shown that the same is true for Gorenstein
rings. Recently, in [{\bf 1}], Bouchiba and Kabbaj showed that if
$R$ and $S$ are $k$--algebras such that $R\otimes_kS$ is
Noetherian then $R\otimes_kS$ is a Cohen--Macaulay ring if and
only if $R$ and $S$ are Cohen--Macaulay rings.

In this paper we shall show that the same is true for complete
intersection and Gorenstein rings. Also it is shown that
$R\otimes_kS$ satisfies Serre's condition $(S_n)$ if and only if
$R$ and $S$ satisfy $(S_n)$.

\vspace{.2in}

\section*{1. Main results}

A Noetherian local ring $R$ is a complete intersection (ring) if
its completion $\hat{R}$ is a residue class ring of a regular
local ring $S$ with respect to an ideal generated by an
$S$--sequence. We say that a Noetherian ring is locally a complete
intersection if all its localizations are complete intersections.

A Noetherian ring $R$ satisfies Serre's condition $(S_n)$ if
$\depth R_{\fp}\ge\Min\{n,\dim R_{\fp}\}$ for all prime ideal
$\fp$ of $R$. Also, a Noetherian ring $R$ satisfies Serre's
normality condition $(R_n)$ if $R_{\fp}$ is a regular local ring
for all prime ideal $\fp$ with $\dim R_{\fp}\le n$.

\vspace{.2in}

The following Theorem is collected from [{\bf 2}; Remark 2.3.5,
Corollary 3.3.15, Theorem 2.1.7, and Theorem 2.2.12]:

\vspace{.1in}

\noindent{\bf Theorem 1.1.} Let $\varphi\,:(R, \fm)\to (S, \fn)$
be a flat local homomorphism of Noetherian local rings. Then the
following hold:

\begin{verse}

(a) $S$ is a complete intersection (resp. Gorenstein,
Cohen--Macaulay) $\Leftrightarrow$ $R$ and $S/\fm S$ are complete
intersections (resp. Gorenstein, Cohen--Macaulay).

(b1) If $S$ is regular then $R$ is regular.

(b2) If $R$ and $S/\fm S$ are regular then $S$ is
regular.\hfill$\square$

\end{verse}

\vspace{.2in}

\noindent{\bf Corollary 1.2.} Let $\varphi\,:R\to S$ be a flat
homomorphism of Noetherian rings. Then the following hold:

\begin{verse}

(a) If $R$ and the fibers $R_{\fp}/\fp R_{\fp}\otimes_RS$,
$\fp\in\Spec(R)$, are regular (resp. locally complete
intersections, Gorenstein, Cohen--Macaulay) then $S$ is regular
(resp. locally complete intersection, Gorenstein,
Cohen--Macaulay).

(b) If $S$ is locally complete intersection (resp. Gorenstein,
Cohen--Macaulay) then the fibres $R_{\fp}/\fp R_{\fp}\otimes_RS$,
$\fp\in\Spec(R)$, are locally complete intersections (resp.
Gorenstein, Cohen--Macaulay)

\end{verse}

\vspace{.1in}

\noindent{\it Proof.} (a): Let $\fq\in\Spec(S)$. Set $\fp=\fq\cap
R\in\Spec(R)$. The induced homomorphism
$\tilde{\varphi}\,:R_{\fp}\to S_{\fq}$ is flat and local. It is
clear that $S_{\fq}/\fp R_{\fp}S_{\fq}$ is a localization of
$R_{\fp}/\fp R_{\fp}\otimes_RS$. Now the assertion follows from
Theorem 1.1.

(b): Let $\fp\in\Spec(R)$. Then $R_{\fp}/\fp
R_{\fp}\otimes_RS\cong S_{\fp}/\fp S_{\fp}$, where
$S_{\fp}=T^{-1}S$ and $T=R-\fp$, and we have

$$\Spec(S_{\fp}/\fp S_{\fp})=\{\fq S_{\fp}/\fp
S_{\fp}|\fq\in\Spec(S), \fq\supseteq\fp S,
\fq\cap(R-\fp)=\varnothing\}.$$

For $\fq S_{\fp}/\fp S_{\fp}\in\Spec(S_{\fp}/\fp S_{\fp})$ we have
to show that $(S_{\fp}/\fp S_{\fp})_{\fq S_{\fp}/\fp S_{\fp}}\cong
S_{\fq}/\fp S_{\fq}$ is complete intersection (resp. Gorenstein,
Cohen--Macaulay). Consider the induced flat local homomorphism
$\tilde{\varphi}\,:R_{\fp}\to S_{\fq}$. Now the assertion follows
from Theorem 1.1.\hfill$\square$

\vspace{.2in}

\noindent{\bf Theorem 1.3.} (See [{\bf 2}; Propositions 2.1.16 and
2.2.21] Let $\varphi\,:R\to S$ be a flat homomorphism of
Noetherian rings. Then the following hold:

\begin{verse}

(a) Let $\fq\in\Spec(S)$ and $\fp=\fq\cap R$. If $S_{\fq}$
satisfies $(S_n)$ (resp. $(R_n)$) then $R_{\fp}$ satisfies $(S_n)$
(resp. $(R_n)$).

(b) If $R$ and the fibers $R_{\fp}/\fp R_{\fp}\otimes_RS$,
$\fp\in\Spec(R)$, satisfy $(S_n)$ (resp. $(R_n)$) then $S$
satisfies $(S_n)$ (resp. $(R_n)$).\hfill$\square$

\end{verse}

\vspace{.2in}

\noindent{\bf Corollary 1.4.} Let $\varphi\,:R\to S$ be a
faithfully flat homomorphism of Noetherian rings. Then the
following hold:

\begin{verse}

(a) If $S$ is regular (resp. locally complete intersection,
Gorenstein, Cohen--Macaulay), then so is $R$.

(b) If $S$ satisfies $(S_n)$ (resp. $(R_n)$), then so does $R$.

\end{verse}

\vspace{.1in}

\noindent{\it Proof.} Let $\fp\in\Spec(R)$. Since $\varphi$ is
faithfully flat there exists $\fq\in\Spec(S)$ such that
$\fp=\fq\cap R$. Consider the flat local homomorphism
$\tilde{\varphi}\,:R_{\fp}\to S_{\fq}$ where
$\tilde{\varphi}(r/s)=\varphi(r)/\varphi(s)$. Now the assertion
follows from Theorems 1.1 and 1.3.\hfill$\square$

\vspace{.2in}

\noindent{\bf Proposition 1.5.} Let $k$ be a field, $L$ and $K$ be
two extension fields of $k$. Suppose that $L\otimes_kK$ is
Noetherian. Then the following hold:

\begin{verse}

(a) $L\otimes_kK$ is locally complete intersection.

(b) If $k$ is perfect then $L\otimes_kK$ is regular.

\end{verse}

\vspace{.1in}

\noindent{\it Proof.} (a): With the same method in the proof of
[{\bf 4}; Theorem 2.2], we can assume that $K$ is a finitely
generated extension field of $k$ (note that, in view of Theorem
1.1, [{\bf 4}; Lemma 2.1] is true with ``Gorenstein ring''
replaced by ``complete intersection''). Now using [{\bf 2};
Proposition 2.1.11] we have that $L\otimes_kK$ is isomorphic to
$$A=T^{-1}(L[x_1,x_2,\ldots ,x_n])/(f_1,f_2,\ldots
,f_m)T^{-1}(L[x_1,x_2,\ldots ,x_n]),$$ where $T$ is a
multiplicatively closed subset of $L[x_1,x_2,\ldots ,x_n]$ and
$f_1,f_2,\ldots ,f_m$ is a $T^{-1}(L[x_1,x_2,\ldots
,x_n])$--sequence. Therefore $A$ is locally complete intersection,
cf. [{\bf 2}; Theorem 2.3.3(c)].

(b): The assertion follows from the note on page 49 of [{\bf 4}]
.\hfill$\square$



\vspace{.2in}

\noindent{\bf Theorem 1.6.} Let $R$ and $S$ be non--zero
$k$--algebras such that $R\otimes_kS$ is Noetherian. Then the
following hold:

\begin{verse}

(a) $R\otimes_kS$ is locally complete intersection (resp.
Gorenstein, Cohen--Macaulay) if and only if $R$ and $S$ are
locally complete intersections (resp. Gorenstein,
Cohen--Macaulay).

(b) $R\otimes_kS$ satisfies $(S_n)$ if and only if $R$ and $S$
satisfy $(S_n)$.

(c) If $R\otimes_kS$ is regular then $R$ and $S$ are regular.

(d) If $R\otimes_kS$ satisfies $(R_n)$ then $R$ and $S$ satisfy
$(R_n)$.

(e) The converse of parts (c) and (d) hold if $\T{char}(k)=0$ or
$\T{char}(k)=p$ such that $k=\{a^p|a\in k\}$.

\end{verse}

\noindent{\it Proof.} Consider two faithfully flat homomorphism:
$$
\varphi\,:R\to R\otimes_kS\,\,\, \mbox{and}\,\,\, \psi\,:S\to
R\otimes_kS$$ \noindent of Noetherian rings.

If $R\otimes_kS$ is regular (resp. locally complete intersection,
Gorenstein, Cohen--Macaulay) then by Corollary 1.4 we have $R$ and
$S$ are regular (resp. locally complete intersections, Gorenstein,
Cohen--Macaulay). Also if $R\otimes_kS$ satisfies $(S_n)$ (resp.
$(R_n)$) then by Corollary 1.4, $R$ and $S$ satisfy $(S_n)$ (resp.
$(R_n)$).

Now let $R$ and $S$ be locally complete intersection (resp.
Gorenstein, Cohen--Macaulay). By Corollary 1.2 it is enough to
show that the fibres $(R\otimes_kS)\otimes_RR_{\fp}/\fp
R_{\fp}\cong R_{\fp}/\fp R_{\fp}\otimes_kS$ over every prime ideal
$\fp$ of $R$ is locally complete intersection (resp. Gorenstein,
Cohen--Macaulay). Consider the flat homomorphism $\gamma\,:S\to
R_{\fp}/\fp R_{\fp}\otimes_kS$. Using Corollary 1.2, it is enough
to show that the fibres $(R_{\fp}/\fp
R_{\fp}\otimes_kS)\otimes_SS_{\fq}/\fq S_{\fq}\cong R_{\fp}/\fp
R_{\fp}\otimes_kS_{\fq}/\fq S_{\fq}$ over every prime $\fq$ of $S$
is locally complete intersection (resp. Gorenstein,
Cohen--Macaulay). But it is clear to see that $R_{\fp}/\fp
R_{\fp}\otimes_kS_{\fq}/\fq S_{\fq}$ is Notherian, since it is a
localization of $R/\fp\otimes_kS/\fq\cong
R\otimes_kS/(\fp\otimes_kS+R\otimes_k\fq)$, which is Noetherian.
Now the assertion follows from Proposition 1.5.

\noindent If $R$ and $S$ satisfy $(S_n)$, with the same proof
$R\otimes_kS$ satisfies $(S_n)$.

\noindent By using the Proposition 1.5 the proof of part (e) is
the same.\hfill$\square$

\vspace{.2in}

\noindent{\bf Remark 1.7.} The converse of part (c) in Theorem 1.6
is not true. For example, let $k$ be an imperfect field of
characteristic 3, let $a\in k$ be an element with no cube root in
$k$. Then $K=k[x]/(x^3-a)k[x]$ is a splitting field of $x^3-a$
over $k$. Thus $K\otimes_kK\cong K[x]/(x^3-a)K[x]$, which is not
regular.


\vspace{.3in}

\noindent {\bf Acknowledgment.} The authors would like to thank
the referee for his/her comments.

\vspace{.3in} \baselineskip=16pt

\begin{center}
\large {\bf References}
\end{center}
\vspace{.2in}

\begin{verse}

[1] Samir Bouchiba and Salah-Eddine Kabbaj, {\em Tensor products
of Cohen-Macaulay rings: solution to a problem of Grothendieck},
J. Algebra {\bf 252} (2002), 65--73.

[2] Winfried Bruns and J\"{u}rgen Herzog, {\em Cohen-Macaulay
rings}, Cambridge University Press, Cambridge, 1993.

[3] A. Grothendieck, {\em \'{E}l\'{e}ments de g\'{e}om\'{e}trie
alg\'{e}brique} IV. \'{E}tude locale des sch\'{e}mas et des
morphismes desch\'{e}mas. II. Inst. Hautes \'{E}tudes Sci. Publ.
Math., No. 24, 1965.

[4] Rodney Y. Sharp, {\em Simplifications in the theory of tensor
products of field extensions.} J. London Math. Soc., {\bf 15},
(1977), 48--50.

[5] Kei-ichi Watanabe, Takeshi Ishikawa, Sadao Tachibana, and Kayo
Otsuka, {\em On tensor products of Gorenstein rings}, J. Math.
Kyoto Univ., {\bf 9} (1969), 413--423.

\end{verse}

\end{document}